\documentstyle[12pt]{article}
\begin{document}
\begin{center}
{\bf \Large The structure of stable
minimal hypersurfaces in $I\!\!R^{n+1}$}
\end{center}
\begin{center}
\begin{tabular}{ccc}
\large Huai-Dong Cao & \large Ying Shen & \large Shunhui Zhu\\
\end{tabular}
\end{center}
\footnotetext[3]{The first author is supported in part by NSF grant
\#DMS-9504925. The third author is supported in part by NSF grant
\#DMS-9404263.}
\begin{abstract}
We provide a new  topological obstruction for
complete stable minimal hypersurfaces
in $I\!\!R^{n+1}$. For $n\geq 3$, we prove that a  complete orientable
 stable minimal hypersurface 
 in $I\!\!R^{n+1}$ cannot have more than one end by showing the existence
of a  bounded harmonic function based on the  Sobolev inequality for minimal
submanifolds \cite{MS}
and by applying the Liouville theorem for harmonic functions due to Schoen-Yau
\cite{SY}.
\end{abstract}

\vspace{1 cm}

\section{Introduction}

This paper is concerned with the structure of complete
stable minimal hypersurfaces in $(n+1)$-dimensional Euclidean space
$I\!\!R^{n+1}$. A complete oriented minimal submanifold $M$
in $I\!\!R^{n+1}$ is called stable if the second variation of
the volume is non-negative on any compact subset of $M$.
The fundamental result along these lines is the
Bernstein Theorem \cite{B} which says that a complete area-minimizing
graph in $I\!\!R^3$ is a plane. Much work has been devoted to trying to
generalize it in the last thirty
years. From the works of Fleming \cite{F}, De Giorgi
\cite{DG}, Almgren
\cite{A} and J. Simons \cite{SJ},  one knows that Bernstein Theorem is
valid for complete area-minimizing graphs in $I\!\!R^{n+1}$ for $n\leq7$.
Counterexamples
to the theorem for $n\geq 8$ were found by Bombieri-De Giorgi-Giusti \cite{BDG}
and later by Lawson \cite{L}. Since then, there have been
attemps to extend the
above Bernstein Theorem by  assuming that the minimal hypersurface
be stable.
The best result in this direction is due to
Fischer-Colbrie-Schoen
 \cite{FS} and do Carmo-Peng \cite{DP}, who proved that if $M^2$ is a
complete,
oriented and immersed stable minimal surface in $I\!\!R^3$, then $M$ is
a plane. 

However, not much is known for the geometric structure of stable minimal
hypersurfaces in $I\!\!R^{n+1}$ when $n\geq 3$. The only known topological obstruction 
for stable minimal
submanifolds we are aware of was a result of 
Schoen-Yau \cite{SY}, which states that if $M^n$ is a complete stable 
hypersurface in a manifold of non-negative curvature and $D$ is a compact 
domain in $M$ with smooth simply connected boundary, then there is no 
non-trivial homomorphism from $\pi _1(D)$  into the fundamental group of 
a compact manifold with nonpositive curvature. We remark
that Schoen-Yau \cite{SY} proved the same result when $M$ is a complete manifold
with non-negative Ricci curvature.
In this paper we provide a new topological obstruction for
complete stable minimal hypersurfaces
in $I\!\!R^{n+1}$. Our main result can be stated as follows:\\

\noindent{\bf Theorem 1.} {\em For any $n\geq 3$, 
if $M^n$ is a complete non-compact
oriented
stable minimal hypersurface in $I\!\!R^{n+1}$, then $M$  has only one end.}\\

To our knowledge, this result is new even for area-minimizing
hypersurfaces $M^n$ in $I\!\!R^{n+1}$ for $n\geq 8$. Note that the stable
condition in Theorem 1 cannot be dropped, since a catenoid, which is unstable, clearly has
two ends. It is also clear that our result differs from Schoen-Yau's theorem.
For example,  our result says that any manifold of
topological type $N^{n-1}\times I\!\!R$ with $N^{n-1}$ compact cannot be a
stable minimal hypersurface in $I\!\!R^{n+1}$. In fact, Theorem 1 can be 
compared with the similar result of Gromoll-Meyer \cite{GM} for complete 
manifolds with positive Ricci curvature.

The proof of our main theorem relies on the Sobolev ineqality for minimal
submanifolds due to Michael and Simon \cite{MS} and the  Liouville
theorem for harmonic maps due to Schoen and Yau \cite{SY}. One
crucial  step in the proof is to show the existence of a non-trivial
bounded harmonic function with finite Dirichlet energy in case the
minimal hypersurface has more than one end. This is done by using the
Sobolev inequality together with a choice of
cut-off functions based on the fact that the minimal submanifold has more
than one end. We remark that our cut-off function actually has noncompact
support. This non-standard choice
of cut-off functions allows us to aviod assumptions such as volume doubling
properties and volume growth conditions of the ends.

Finally, we would like to point out that the method we used to in the
proof of Lemma 2 in next section yields the following result which is
of independent interest:\\

\noindent{\bf Theorem 2.} {\em Let $M^n$ be a complete noncompact
Riemannian manifold with at least two ends of infinite volume. Suppose that
either

(1) the Sobolev inequality holds on $M$, or 

(2) the first eigenvalue $\lambda_1(M)$ of $M$ is positive.
}\\
{\em Then there exists on $M$ a non-constant bounded harmonic function 
with finite Dirichlet energy. }\\

The problem of the existence of harmonic functions on 
a complete manifold has a long history.  For some of the recent 
progress on this problem, we refer the readers to the 
important works of Li-Tam (\cite{LT1} and \cite{LT2}), Colding-Minicozzi
(\cite{CM1}, \cite{CM2}, \cite{CM3} and \cite{CM4}) and a very 
recent survey article of Peter Li \cite{Li}. 

{\bf Acknowledgement} We would like to thank Prof. Richard 
Schoen for bringing this problem to our attention and for very helpful
discussions. We would also like to  thank  Professors Jeff Cheeger, Peter Li, and S.-T. Yau for their
interest in this work. 

\section{The Proof}

In this section we prove Theorem 1 stated in the introduction. First let us
fix some notations and recall the definition of a stable
minimal submanifold.

Let $\{U_i\}_{i=1}^{\infty}$ be a family of relatively compact open sets
which  exhaust the manifold $M$, i.e.,
\[U_i\subset U_{i+1},\]
\[ \cup_{i=1}^{\infty}U_i=M. \]

An end of $M$ is an inverse system $E=\{E^{(i)}\}_{i=1}^{\infty}$ such that
\[E^{(i+1)}\subset E^{(i)}\]
and $E^{(i)}$ is a connected component of $M\setminus \overline{U}_i$.\\

{\bf Remark 2}. If $M$ has only finitely many ends, then there is
some $i_0 >0$ such that all inverse systems ($E^{(i)}$) stabilize for $i\geq
i_0$,
 i.e., $E^{(i)}=E^{(i_0)}$.\\

{\bf Remark 3}. An end can also be defined as equivalent classes
of cofinal curves, where two curves $c_1, c_2$: $[0, \infty)\rightarrow M$
are cofinal iff for every compact set $K \subset M$ there is some $t>0$
such that $c_1(t_1)$ and $c_2(t_2)$ lie in the same connected component
of $M\setminus K$ for all $t_1, t_2 > t$.\\

Let $M^n$ be a completed oriented submanifold minimally immersed in
$I\!\!R^{n+1}$. We say that $M$ is stable if the second variation of
the volume is non-negative on any compact subset of $M$.
More precisely,
let $e_1, e_2, \cdots, e_{n+1}$ be a positively oriented orthonormal
frame on $M$ with $e_1, \cdots, e_n$ tangential, and $e_{n+1}$ the
globally defined unit normal vector to $M$. We can define the second
fundamental form $\{h_{ij}\}$ of $M$ by
\[      h_{ij}=\langle \nabla_{e_i}e_{n+1}, e_j \rangle \]
for $i,j=1, 2, \cdots, n$, where $\nabla$ is the Riemannian connection of
$I\!\!R^{n+1}$. Then $M$ is minimal iff the mean curvature
\[      H=\sum_{i=1}^{n}h_{ii}=0.\]
The stability of $M$ is given by the following inequality
(see \cite{SJ} or \cite{SY}):
\[      \int_{M} |\nabla \phi|^2 \geq \int_{M} \sum _{i,j=1}^{n} h_{ij}^2
\phi^2 \]
where $\phi$ is any function with compact support on $M$.\\

Now we are ready to present the proof of Theorem 1.\\

\noindent{\bf Lemma 1.} {\em If $M^n$ is a complete orientable minimal
hypersurface in $I\!\! R^{n+1}$, then every end of $M^n$ has infinite
volume.}\\

{\bf Proof}: In fact we will show that for any compact set $K\subset M$, 
every noncompact component of $M\setminus K$ has infinite volume. Let $E$ be
 a component of $M\setminus K$. We will adopt an argument of Yau (\cite{Y})
 to the end $E$.

Take an arbitrary point $p\in M$, without loss of generality, we may assume
$p=0$. In the following we let $d(\cdot, \cdot)$ be the distance function of
$I\!\!R^{n+1}$, and $r(\cdot, \cdot)$ the distance function of $M$ with 
respect 
to the induced metric. We will write $d(x), r(x)$ if the base point is $0$. 
Obviously
$d\leq r$ for any two points in $M$. Let $\gamma$ be a minimal geodesic from
$0$, then,
\begin{eqnarray}
{\partial d \over \partial r} & = &\lim _{t\rightarrow 0}{d(\gamma
(s+t))-d(\gamma (s))\over t} \nonumber\\
& \leq & \lim _{t\rightarrow 0}{d(\gamma
(s+t),\gamma (s))\over t} \;\;\;\;\;(\mbox{ by the triangle inequality})
\nonumber\\ & \leq & 1.\;\;\;\;\;\;\; (\mbox{since}\;\;d \leq r)
\label{1}
\end{eqnarray}

By a direct computation and using the fact that $M$ is minimal, one can show that 
\begin{eqnarray*}
\triangle _M d^2 (x) &= &  2n.
\end{eqnarray*}

Let $B(s)$ be the geodeisc ball of $M$,  of radius $s$ centered at $0$.
Integrating 
the above equation over $B(s)$ and using (\ref{1}), we obtain

\begin{eqnarray*}
2n\; \mbox{vol}(B(s)) 
&\leq & 2s\; \mbox{vol}(\partial B(s)).
\end{eqnarray*}

Note that in any manifold, 
\[\mbox{vol}(\partial B(s)) ={\partial\over \partial r} \mid _{r=s}
\mbox{vol}(B(s)) .\] We thus obtain
\[s\;{\partial \over \partial r}|_{r=s} \mbox{vol}(B(r))- n\;
\mbox{vol}(B(s))\geq 0,\]
which implies $s^{-n}\mbox{vol}(B(s))$ is nondecreasing. Therefore
\[{\mbox{vol}(B(s)) \over s^n}\geq \lim _{s\rightarrow 0}{\mbox{vol}(B(s))
\over s^n}
= \omega (n)\]
where $\omega (n)$ is the volume of unit ball in $I\!\! R^{n}$.

Now if $E$ has finite volume, choose $R$ big enough such that
\[\omega (n) R^n> \mbox{vol}(E).\]
Let $p$ be a point in $E$ such that $r(p, \partial E)\geq R$, then
\[\mbox{vol}(E) \geq \mbox{vol}(B(R)) \geq \omega (n) R^n> \mbox{vol}(E),\]
a contradiction.
\hspace*{\fill}q.e.d\\

\noindent{\bf Lemma 2.} {\em Let $M^n$ be a complete orientable minimal
hypersurface in $I\!\! R^{n+1}$ with at least two ends. 
Then there exists on $M$ a non-constant bounded harmonic function
with finite energy.}\\

We remark that the statement of Lemma 2. still holds when $M^n$ is a complete submanifold in
$I\!\! R^{n+m}$ for $n\geq 3$, $m\geq 1$ provided one of the following 
conditions holds:
(1) $\int_M |\vec {H}|^n < 1/c_n^{\frac {n}{2}}$; or 
(2) $\int_M |\vec {H}|^{\frac {4n}{n+2}} < \infty$.
Here $\vec{H}$ is the mean curvature vector
of $M$ and $c_n$ is a (Sobolev) constant which depends only on the dimension
$n$. 

{\bf Proof of Lemma 2}: In the following, we will use an exhausion of $M$ by
compact submanifolds with boundary.  The obvious choice is the distance
balls which are not smooth in general. It is standard to smooth the
distance function to a point $p$  by, say, local averaging,
to get a smooth function $f$.
$f$ still has compact sublevel sets. By Sard's Theorem, we can choose a
sequence of regular values $\{R_i\}$ of $f$, such that $\{ f^{-1}(0,
R_i)\}$ gives an exhausion of $M$ by
compact smooth submanifolds. In the following, we denote $D_i=f^{-1}(0,
R_i)$. We shall also use the notations for the ends of $M$ introduced at
the begining of this
section. For $i\geq i_0$ and $i_0$ sufficiently large, let
\[M\setminus D_i=\cup_{j=1}^{s} E_j^{(i)}\]
be the disjoint union of connected components, with $s\geq 2$. 

By Lemma 1 and the assumptions in Lemma 2, $M$ has at least two components with infinite 
volume. Let $E_1^{(i_0)},
E_2^{(i_0)}$ be two such components. On each compact domain $D_i$, we can
minimize
the engergy functional
$\int_{D_i} |\nabla u|^2 dx$
among all functions
$u$ such
that $u|_{\partial E_1^{(i)}}=1$ and $u|_{\partial E_j^{(i)}}=0$ for
all $j\geq 2$, where $\nabla$ and $dx$ are gradient and volume element of
$M$, respectively.  We denote the minimizer by $u_i$. Then  $u_i$ is
the unique solution of the following Dirichlet problem on $D_i$
\begin{equation}
\left\{
\begin{array}{lcl}
   \triangle u(x) &= &0\\
u|_{\partial E_1^{(i)}}&= &1  \\
u|_{\partial E_j^{(i)}}&= &0.\;\; (j\neq 1)
\label{3}
\end{array}\right.
\end{equation}

By maximum principle, we have $0\leq u_i\leq 1$ on $D_i$. Moreover, it is
easy to see that $\int_{D_i} |\nabla u_i|^2 dx
\leq \int_{D_j} |\nabla u_j|^2 dx$ for $i>j$ and hence
there is a universal constant $C_1>0$ such that
\begin{eqnarray}
\int_{D_i} |\nabla u_i|^2 dx < C_1.
\label{4}
\end{eqnarray}

Therefore by passing to a subsequence, still denoted by $u_i$, we can
find a harmonic function $u$ on $M$ such that
\[	\lim_{i\rightarrow \infty} u_i(x) = u(x), \hspace{1.0 cm} x\in M \]

 and \[	\int_M |\nabla u|^2 dx <C_1.  \] It is clear from the construction
that $0\leq u(x) \leq 1$ on $M$.

In the following we prove that the limiting harmonic
function $u$ is not a
constant function.
We will prove this by contradiction.

By using the Sobolev inequality for minimal hyperfaces in $I\!\! R^{n+1}$ (see
\cite {MS}),
we have, for any smooth function $\phi $ which vanishes on $\partial D_i$,
\begin{eqnarray}
(\int_{D_i}  \phi^p(x)dx)^{\frac {2}{p}} \leq c_n \int_{D_i}|\nabla 
\phi|^2 dx 
\label{5}
\end{eqnarray}
where $p=2n/n-2$ (this is where we need to assume $n\geq 3$) and
$c_n$ is the Sobolev constant which only depends on the dimesion $n$.

Note that from the construction of $u_i$, the function $u_i
(1-u_i)$ vanishes on $\partial D_i$. Setting $\phi =u_i(1-u_i)$ in
(\ref{5}),  we  obtain
\begin{eqnarray}
(\int_{D_i} ( u_i (1-u_i) )^p dx)^{\frac {2}{p}}&\leq &
c_n\int _{D_i} (\nabla u_i-2u_i\nabla u_i)^2 dx \nonumber\\
&\leq & c_n\int _{D_i} (|\nabla u_i|+2u_i|\nabla u_i|)^2 dx \nonumber\\
&\leq & 9 c_n\int_{D_i}|\nabla u_i|^2 dx \leq 9c_n C_1.
\label{6}
\end{eqnarray}

Since $\mbox {vol}(D_i)\rightarrow \infty$, by letting $i\rightarrow
\infty$ in (\ref{6}), it follows that
if $u$  is
a constant function,  then $u\equiv 0$ or $u\equiv 1$.

Thus we only need to show that $u$ cannot be identically $0$ or
$1$.
If $u\equiv 0$, we may replace $u_i$ and $u$ by
$\hat{u}_i=1-u_i$ and $\hat{u}=1-
u$, respectively. Then $\hat{u}\equiv 1$ and furthermore, $\hat{u}_i$ satisfies
 (\ref {4}), and (\ref {6}). Thus we  may assume that
$u\equiv 1$. (This is where we use the condition that the manifold has
 at least two ends with infinite volume.)

We choose a smooth function $\psi$ such that
\[
\psi = \left\{
\begin{array}{ll}
1 & \mbox{in}\;\; E_2^{(i_0)} \\
0 & \mbox{in}\;\; E_j^{(i_0)}\;\;j\neq 2
\end{array}
\right. \]
and
\[	 |\nabla \psi|\leq C_2 ,\;\;\; 0\leq \psi \leq 1 \]
for some constant $C_2>0$ which is independent of $i$ and $u_i$. Note 
that $ |\nabla \psi|$ vanishes outside a compact set.

Since $u_i|_{ E_1^{(i)}}= 1$ and $u_i|_{ E_j^{(i)}}= 0$
for $j\geq 2$, the function $\phi =u_i\psi$ vanishes on $\partial D_i$.
The Sobolev inequality (5) implies that
\begin{eqnarray}
(\int_{D_i} (u_i \psi)^p dx)^{\frac {2}{p}}& \leq & c_n\int_{D_i} |\nabla
(u_i\psi)|^2 dx   \nonumber \\
& \leq & c_n\int_{D_i} |\psi \nabla
u_i|+u_i |\nabla \psi |)^2 dx \nonumber \\
& \leq & 2 c_n ( \int_{D_i}\psi^2 |\nabla u_i|^2 +
\int_{D_i} |\nabla \psi|^2) \leq C_3
\label{7}
\end{eqnarray}
where $C_3$ is a constant independent of $i$ and $u_i$.

Therefore, we have

\begin{eqnarray*}
        \int_{E_2^{(i_0)}\cap D_i} (u_i)^p &\leq C_3
\end{eqnarray*}
for all $i\geq i_0$.
Letting $i\rightarrow \infty$, we get
\begin{eqnarray}
\mbox{vol}(E_2^{(i_0)})=
	\int_{E_2^{(i_0)}} u^p\leq C_3.
\end{eqnarray}

This contradicts our assumption that $\mbox{vol}(E_2^{(i_0)})$ is
infinite. Therefore the limiting harmonic function $u$ is not a constant 
function.
\hspace*{\fill}q.e.d. \\

\noindent{\bf Lemma 3.} (Schoen-Yau \cite{SY}) {\em Let $M$ be a complete
noncompact
stable minimal hypersurface in a manifold of non-negative curvature. If
$u$ is a harmonic function on $M$ with bounded energy, then $u$ is constant.}\\

This is a special case of the Liouville Theorem for harmonic map that
Schoen and Yau originally proved. 

 Now suppose $M^n (n\geq 3)$ is a complete, oriented
stable minimal hypersurface in $I\!\!R^{n+1}$ with finitely many ends.
From Lemma 1 we know that each end of $M$ has infinite volume. If $M$ 
has more than one end then Lemma 2 implies that $M$ supports a 
non-constant harmonic function with
finite energy. This is a contradiction to  Lemma 3. Hence $M$ has only one 
end and the proof of Theorem 1 is completed.

\noindent Department of Mathematics, Texas A$\&$M University, College Station,
TX
77843\\
{\em E-mail address}: cao@math.tamu.edu\\

\noindent Department of Mathematics, Dartmouth College, Hanover, NH 03755\\
{\em E-mail address}: Ying.Shen@dartvax.dartmouth.edu\\

\noindent Department of Mathematics, Dartmouth College, Hanover, NH 03755\\
{\em E-mail address}: Shunhui.Zhu@dartvax.dartmouth.edu

\end{document}